\documentclass[12pt]{amsart}

\usepackage{mathrsfs}
\usepackage{amsmath,amsthm,amssymb}
\font\teneufm=eufm10
\font\seveneufm=eufm7
\font\fiveeufm=eufm5
\newfam\frakturfam

\textfont\frakturfam=\teneufm
\scriptfont\frakturfam=\seveneufm
\scriptscriptfont\frakturfam=\fiveeufm

\newtheorem{pr}{Proposition}
\newtheorem{df}{Definition}
\newtheorem{lm}{Lemma}
\newtheorem{theor}{Theorem}
\newtheorem{co}{Corollary}

\newtheorem{prob}{Problem}
\def\bee{\begin{eqnarray}}
\def\bes{\begin{eqnarray*}}
\def\eee{\end{eqnarray}}
\def\ees{\end{eqnarray*}}

\def\a{\alpha}
\def\b{\beta}

\def\s{\sigma}

\def\d{\partial}
\def\Proof{{\sl Proof.}\ }

\title{Tangent Lie Algebras of Automorphism Groups of Free Algebras}

\begin{document}
\date{}
\maketitle

\begin{center}
{\bf Ivan Shestakov}\footnote{Instituto de Matem\'atica e Estat\'\i stica, Universidade de S\~ao Paulo, 
Caixa Postal 66281, S\~ao Paulo - SP, \mbox{05315--970}, Brazil 
e-mail:{\em shestak@ime.usp.br}}
and
{\bf Ualbai Umirbaev}\footnote{Department of Mathematics,
 Wayne State University,
Detroit, MI 48202, USA 
and Institute of Mathematics and Mathematical Modeling, Almaty, 050010, Kazakhstan
e-mail: {\em umirbaev@wayne.edu}}
\end{center}

\begin{abstract} 
We study an analogue of the Andreadakis–Johnson filtration for automorphism groups of free algebras and introduce the notion of tangent Lie algebras for certain automorphism groups, defined as subalgebras of the Lie algebra of derivations. We show that, for many classical varieties of algebras, the tangent Lie algebra is contained in the Lie algebra of derivations with constant divergence. 
We also introduce the concepts of approximately tame and absolutely wild automorphisms of free algebras in arbitrary varieties and employ tangent Lie algebras to investigate their properties. It is shown that nearly all known examples of wild automorphisms of free algebras are absolutely wild, with the notable exceptions of the Nagata and Anick automorphisms. We show that the Bergman automorphism of free matrix algebras of order two is absolutely wild. Furthermore, we prove that free algebras in any variety of polynilpotent Lie algebras--except for the abelian and metabelian varieties--also possess absolutely wild automorphisms.
\end{abstract}



\tableofcontents

\section{Introduction}

\hspace*{\parindent}

Let $F_n$ be the free group of rank $n$ and let $H=F_n^{\mathrm{ab}}=F_n/[F_n,F_n]$ be its abelianization. The natural homomorphism  $F_n\to H$ induces a homomorphism  
\bes
\rho : \mathrm{Aut}(F_n)\to \mathrm{Aut}(H)=GL_n(\mathbb{Z})
\ees
between the corresponding automorphism groups. The elements of the kernel  $\mathrm{IA}_n=\mathrm{Ker}(\rho)$ of this homomorphism are called  {\em $\mathrm{IA}$ automorphisms} of $F_n$. The group $\mathrm{IA}_n$ forms a major part of $\mathrm{Aut}(F_n)$, and its structure is very complex  \cite{Bestvina23,Vogtmann18}.

Consider the lower central series  
\bes
F_n=\Gamma_n(1)\supseteq\Gamma_n(2)\supseteq\ldots\supseteq \Gamma_n(k)\supseteq
\ees
of   $F_n$, where $\Gamma_n(i+1)=[\Gamma_n(i),\Gamma_n(1)]$ for all $i\geq 1$. The action of $\mathrm{Aut}(F_n)$ on the quotient $F_n/\Gamma_n(k+1)$ is well defined,  since $\Gamma_n(k)$ is invariant under the action of $\mathrm{Aut}(F_n)$. The natural homomorphism $F_n\to F_n/\Gamma_n(k+1)$ induces a homomorphism  
\bes
\rho_k : \mathrm{Aut}(F_n)\to \mathrm{Aut}(F_n)/\Gamma_n(k+1), \ \ k\geq 1. 
\ees
Let  $\mathcal{A}_n(k)=\mathrm{Ker}(\rho_k)$. The descending central series   
\bes
\mathrm{IA}_n=\mathcal{A}_n(1)\supseteq\mathcal{A}_n(2)\supseteq\ldots\supseteq \mathcal{A}_n(k)\supseteq
\ees
is called the  {\em Andreadakis--Johnson filtration} of $\mathrm{Aut}(F_n)$ \cite{Satoh16}. 

 The graded quotients $\mathrm{gr}^k( \mathcal{A}_n)= \mathcal{A}_n(k)/ \mathcal{A}_n(k+1)$ naturally have  $GL_n(\mathbb{Z})$-module structures and the direct sum  
\bes
\mathrm{gr}( \mathcal{A}_n)=\mathrm{gr}^1( \mathcal{A}_n)\oplus \mathrm{gr}^2( \mathcal{A}_n)\oplus\ldots\oplus \mathrm{gr}^k( \mathcal{A}_n)\oplus\ldots
\ees
 naturally can be turned into a Lie ring. Its elements can be interpreted as derivations of a free Lie ring via the so-called Johnson homomorphisms, and many interesting papers are devoted to studying the structure of this ring \cite{Satoh16}. 

These studies are not very popular among specialists in ring theory. Nevertheless, similar constructions have also been studied by specialists in ring theory: by Anick \cite{Anick} in the context of polynomial algebras over a field of characteristic zero, and by Bryant and Drensky \cite{BD93-1} for free metabelian Lie algebras over a field of characteristic zero. However, they did not connect these Lie algebras to derivations of free algebras.

 Let ${\mathfrak M}$ be an arbitrary variety of algebras over a filed $K$ of characteristic zero and let 
 $A=K_{\mathfrak M}\langle x_1,x_2,\ldots,x_n\rangle$ be the free algebra of ${\mathfrak M}$ with a free set of generators $X=\{x_1,x_2,\ldots,x_n\}$. Let $\mathrm{Aut}(A)$ be the group of all automorphisms of $A$.

In this paper, we study the subgroups $H$ of $\mathrm{Aut}(A)$ that contain the subgroup  $\mathrm{Aff}_n(K)$ of all affine automorphisms when $\mathfrak{M}$ is a unital variety of algebras and $A$ is a unital algebra, and that contain only the subgroup of linear automorphisms $GL_n(K)$ when $A$ is not unital. This concerns the general approach to describing generators of certain automorphism groups modulo  $\mathrm{Aff}_n(K)$ or $GL_n(K)$. Some results of this type are given in Section 4.  

For these type automorphism groups $H$ we define the tangent Lie algebra $T(H)$ with respect to the power series topology. 
The tangent Lie algebras are directly defined as subalgebras of the Lie algebra $\mathrm{Der}(A)$ of all derivations of $A$. The tangent Lie algebra $T(\mathrm{Aut}(A))$ is a slight extension of the algebra corresponding to $\mathrm{gr}( \mathcal{A}_n)$, but it is more convenient for formulating certain problems concerning automorphism groups. 

We use tangent algebras to address the problem of characterizing tame and wild automorphisms. Denote by $\phi = (f_1,f_2,\ldots,f_n)$ the automorphism 
  $\phi$ of $A$ such that $\phi(x_i)=f_i,\, 1\leq i\leq n$. An automorphism  
\bes
(x_1,\ldots,x_{i-1}, \a x_i+f, x_{i+1},\ldots,x_n),
\ees
where $0\neq\a\in K,\ f\in K_{\mathfrak M}\langle X\setminus \{x_i\}\rangle$, 
 is called {\em elementary}.    The subgroup $\mathrm{TAut}(A)$ of $\mathrm{Aut}(A)$ generated by all 
  elementary automorphisms is called the {\em tame automorphism group}, 
 and the elements of this subgroup are called {\em tame automorphisms} 
 of $A$. Nontame automorphisms of $A$ are called {\em wild}.

The well known  Jung--van der Kulk Theorem \cite{Jung, Kulk} says that all automorphisms of the polynomial algebra $K[x,y]$ in two variables $x,y$ over a field $K$ are tame. Similar results hold for free associative algebras \cite{Czer, ML70} and for free Poisson algebras (in characteristic zero) \cite{MLTU} (see also \cite{MLSh,U12,MLU16}).  Moreover, the automorphism groups of polynomial algebras, free associative algebras, and free Poisson algebras in two variables (over a field of characteristic zero) are isomorphic.

The automorphism  groups of commutative and associative algebras generated by three elements are much more complicated. 
The well-known Nagata automorphism (see \cite{Nagata})
\bes
(x+2y(zx-y^2)+z(zx-y^2)^2, y+z(zx-y^2),z)
\ees
of the polynomial algebra $K[x,y,z]$ over a field $K$ of characteristic $0$ is proven to be wild \cite{SU04-1,US02RAN}. 
The Anick automorphism (see \cite[p. 398]{Cohn06})
\bes
(x+z(xz-zy),\,y+(xz-zy)z,\,z)
\ees
of the free associative algebra $K\langle x,y,z\rangle$ over a field $K$ of characteristic $0$ is also proven to be wild \cite{U07JR,06AN}. 
The Nagata
automorphism gives an example of a wild
automorphism of free Poisson algebras in three variables. Recently Shestakov and Zhang \cite{SZ24} constructed an analogue of the Anick automorphism for free Poisson algebras in three variables. 

It is well known \cite{Smith} that the Nagata and Anick automorphisms are stably tame, that is, they become tame after adding one more variable. 

 In 1964 P. Cohn  proved \cite{Cohn64} that all automorphisms of finitely generated free
Lie algebras over a field are tame.
 Later this result was extended to free algebras of Nielsen-Schreier varieties \cite{Lewin}. Recall that
a variety of universal algebras is called Nielsen-Schreier, if any subalgebra of a free algebra of this variety is free,
i.e., an analog of the classical Nielsen-Schreier theorem is true.
 The varieties of all non-associative algebras \cite{Kurosh},
 commutative and anti-commutative algebras \cite{Shirshov54}, Lie algebras and Lie $p$-algebras 
 \cite{Shirshov53,Witt}, and Lie superalgebras  \cite{Mikhalev85,Stern} and Lie $p$-superalgebras \cite{Mikhalev88} over a field are Nielsen-Schreier. 
 Some other examples of Nielsen-Schreier 
 varieties can be found in \cite{SU02,MS14,U94,U96,Chibrikov}. 
 
It was recently shown \cite{DU22} that the varieties of pre-Lie (also known as right-symmetric) algebras and Lie-admissible algebras over a field of characteristic zero are Nielsen--Schreier. In particular, every automorphism of a free right-symmetric and a free Lie-admissable algebra of finite rank is tame. Tameness of automorphisms of free right-symmetric algebras in two variables over a field any characteristic was proven in \cite{KMLU}. 

Defining relations for automorphism groups of finitely generated free algebras of Nielsen-Schreier varieties were described in \cite{U07JA}.

An automorphism $\phi$ will be called {\em approximately tame} if it can be approximated by a sequence of 
tame automorphisms $\{\psi_k\}_{k\geq 0}$ with respect to the {\em formal power series topology}  (see exact definitions in Section 4).  An automorphism $\phi$ will be called {\em absolutely wild} if it is not approximately tame.

In 1983 Anick \cite{Anick} proved that every automorphism of the polynomial algebra $K[x_1,x_2,\ldots,x_n]$ over a field $K$ of characteristic zero is approximately tame. 
He  considered the power series topology and applied it for endomorphisms too. Anick considered his result as a weak generalization of the above mentioned Jung--van der Kulk Theorem.

This means that the Nagata automorphism is wild but approximately tame.
 An analogue of the Anick result for free associative algebras and free Poisson algebras is still unknown. In particular, we don't know weather if the above mentioned wild Anick automorphism and its analogue for Poisson algebras are approximately tame or not.

\begin{prob}\label{pr2} Is every automorphism of a free associative algebra and a free Poisson algebra in $n\geq 3$ variables over a field of characteristic zero approximately tame?
\end{prob}

We show that the tangent Lie algebras in many cases coinside with the {\em special Lie algebras of derivations}, that is, the algebra of derivations without divergence. 
The notion of divergence of derivations was known in the cases of free Lie algebras (see, for example, Enomoto and Satoh \cite{ES11} and Satoh \cite{Satoh12}) and in the case of free associative algebras (see, for example, Alekseev, Kawazumi, Kuno and Naef \cite{AKKN}). Recently, the notion of divergence for a derivation of an arbitrary free operadic algebra was defined by Powell \cite{Powell21}; his generalized divergence takes values in the commutator quotient of the universal multiplicative enveloping algebra. 

We demonstrate how the study of approximately tame and absolutely wild automorphisms can be related to the study of generators of the tangent Lie algebras. Using this we reformulate the results of Bryant and Drensky \cite{BD93-1} and by Kofinas and Papistas \cite{KP16}. Note that all examples of wild automorphisms given in \cite{BSh95,BD93-1,Papistas93} are absolutely wild. In the language of tangent algebras those wild automorphisms can be detected algorithmically.  

This paper is organized as follows. In Section 2, we define tangent Lie algebras and describe their elements using the language of tangent derivations. In Section 3, employing universal derivations, Jacobian matrices, and divergence, we show that the tangent Lie algebra, in certain important cases, is a subalgebra of the Lie algebra of derivations with constant divergence. In Section 4, we define approximately tame and absolutely wild automorphisms and reformulate several results on automorphisms in terms of tangent algebras. Finally, in Section 5, we demonstrate some methods for detecting absolutely wild automorphisms.

\section{Tangent Lie algebras}

\hspace*{\parindent}

Let ${\mathfrak M}$ be an arbitrary variety of algebras over a filed $K$ of characteristic zero and let 
 $A=K_{\mathfrak M}\langle x_1,x_2,\ldots,x_n\rangle$ be the free algebra of ${\mathfrak M}$ with a free set of generators $X=\{x_1,x_2,\ldots,x_n\}$. If ${\mathfrak M}$ is a unitary variety of algebras (see, for example \cite{KBKA}) then we suppose that $A$ has the identity element $1$. 
Consider the natural grading 
\bee\label{f1}
 A=A_0\oplus A_1\oplus \ldots \oplus A_k\oplus\ldots 
\eee
with respect to the standard function $\deg$, that is, $A_k$ is the linear span of (non-associative) monomials of degree $k$. An non-zero element $f\in A_k$ is called {\em homogeneous} of degree $k$. We have $A_0=K\cdot 1$ if $A$ is unital, and $A_0=0$ otherwise. 

Let $L=\mathrm{Der} (A)$ be the Lie algebra of all derivations of $A$. For any $n$-tuple $F=(f_1,\ldots,f_n)\in A^n$ there exists a unique derivation $D$ of $A$ such that $D(x_i)=f_i$ for all $1\leq i\leq n$. Denote this derivation by 
\bee\label{f4}
D=f_1\partial_1+\ldots+f_n\partial_n=D_F. 
\eee
We say that $D$ is homogeneous of degree $i$ if $f_1,\ldots,f_n\in A_{i+1}$. Let $L_i$ be the space of all homogeneous derivations of degree $i$. Then 
\bes
L=L_{-1}\oplus L_0\oplus L_1\oplus\ldots\oplus L_k\oplus \ldots, \ \ \ [L_i,L_j]\subseteq L_{i+j}, 
\ees
is a grading of $L$. We have $L_{-1}=K\partial_1+\ldots+K\partial_n$ if $A$ has an identity element and $L_{-1}=0$ otherwise, and $L_0$ is isomorphic to the matrix algebra $\mathfrak{gl}(n)$ with the matrix units $e_{ij}=x_i\partial_j$, where $1\leq i,j\leq n$. 

We also consider homogeneous with respect to one variable $x_i$ elements of $A$ and derivations. Recall that if $f\in A$ is homogeneous with respect to $x_i$ of degree $\mathrm{deg}_{x_i}(f)=k$ then $f\d_j$ is homogeneous of degree with respect to $x_i$ of degree $\mathrm{deg}_{x_i}(f\d_j)=k-\delta_{ij}$, where $\delta_{ij}$ is the Kronecker delta. 

The {\em left-symmetric} product 
\bes
u\d_i\cdot v\d_j=(u\d_i)(v)\d_j 
\ees
turns the space of all derivations $L$ into a left symmetric algebra \cite{CA17}. Then $[D_1,D_2]=D_1\cdot D_2-D_2\cdot D_1$ for any $D_1,D_2\in D$.

We are going to study some subalgebras of $L$ containing $L_0$. The following proposition describes some important properties of such subalgebras. Recall that every automorphism of $A$ acts by conjugation on the algebra of derivations $L$.

\begin{pr}\label{X} Let $R$ be a subalgebra of $L$ containing $L_0$. Then 

$(a)$ $R$  is homogeneous with respect to each variable $x_i$. In particular, $R$ is homogeneous with respect to the degree function $\mathrm{deg}$. 

$(b)$ $R$ is invariant under the action of the group of linear automorphisms $GL_n(K)$. 
\end{pr}
\Proof $(a)$ Let $D\in R$ and let 
\bes
D=D_{-1}+D_0+D_1+\ldots+D_k
\ees
be the homogeneous decomposition of $D$ with respect $x_1$ with  $\deg_{x_1}(D_i)=i$ for all $-1\leq i\leq k$. 
We have $[x_1\d_1,D_i]=iD_i$ for all $i$. Consequently, 
\bes
\mathrm{ad}(x_1\d_1)^q(D)=\sum_{i=-1}^k i^qD_i\in R
\ees
for all $q\geq 0$. An application of the Vandermonde determinant argument gives that $D_i\in R$ for all $i$. 

$(b)$ It is well known that $GL_n(K)$ is generated by all automorphisms of the form $(\a x_1,x_2,\ldots,x_n)$, $0\neq\a\in K$,  and all elementary automorphisms $(x_1,\ldots,x_i+t x_j,\ldots,x_n)$, $t\in K$. 

First set $\phi=(\a x_1,x_2,\ldots,x_n)$. We want to show that $\phi D \phi^{-1}\in R$ for any $D\in R$. By $(a)$, we may assume that $D$ is homogeneous with respect to $x_1$ of degree $\mathrm{deg}_{x_1}(D)=k$. Then $\phi D \phi^{-1}=\a^kD\in R$.

Let $\phi=(x_1+t x_2,x_2,\ldots,x_n)$ and let $D=f_1\d_1+f_2\d_2+\ldots+f_n\d_n\in R$ be a homogeneous with respect to $x_1$ derivation of degree $\mathrm{deg}_{x_1}(D)=k$. We have  
\bes
\phi D\phi^{-1}= (\phi(f_1)-t\phi(f_2))\d_1+\phi(f_2)\d_2+\ldots+\phi(f_n)\d_n. 
\ees
Set $s=x_2\d_1$. For any $f\in A_{k+1}$ we have 
\bes
\phi(f)=f+ts(f)+\frac{1}{2!}t^2s^2(f)+\ldots+\frac{1}{(k+1)!}t^{+1}s^{k+1}(f). 
\ees
Consequently, 
\bes
\phi D\phi^{-1}=D+tsD+\frac{1}{2!}t^2s^2D+\ldots+\frac{1}{(k+1)!}t^{k+1}s^{k+1}D\\
-t(f_2+ts(f_2)+\frac{1}{2!}t^2s^2(f_2)+\ldots+\frac{1}{(k+1)!}t^{k+1}s^{k+1}(f_2))\d_1, 
\ees
where $s^{i+1}D=s\cdot s^iD$ for all $i\geq 0$ is defined inductively using the left-symmetric product introduced above. Then
\bes
\phi D\phi^{-1}=D+t(sD-f_2\d_1)+t^2(\frac{1}{2!}s^2D-s(f_2)\d_1)+\ldots\\
+t^{k+1}(\frac{1}{(k+1)!}s^{k+1}D-\frac{1}{k!}s^{k}(f_2)\d_1)-\frac{1}{(k+1)!}t^{+1}s^{k+1}(f_2)\d_1. 
\ees
It is easy to verify by induction on $q$ that
\bes
\mathrm{ad}(s)^q(D)=s^qD-qs^{q-1}(f_2)\d_1, q\geq 0. 
\ees
Consequently, 
\bes
\phi D\phi^{-1}=\sum_{p=0}^{k+2}\frac{1}{p!} \mathrm{ad}(s)^pD\in R
\ees
since $s^{k+2}D=0$. $\Box$

Let $\mathrm{Aut}(A)$ be the group of all automorphisms of $A$. Let $\mathrm{IA}(k)=\mathrm{IA}(A,k)$ be the 
 set of all automorphisms of $A$ that induces the identity automorphism on the factor-algebra $A/(A_{k+1}+A_{k+2}+\ldots)$. 
The group of automorphisms $\mathrm{IA}(1)=\mathrm{IA}(A)$ is called the group of {\em $\mathrm{IA}$-automorphisms} of $A$.   
 
Consider the descending central series 
\bee\label{f2}
\mathrm{IA}(A)=\mathrm{IA}(1)\supseteq\mathrm{IA}(2)\supseteq\ldots\supseteq \mathrm{IA}(k)\supseteq . 
\eee
 
Let $\phi\in \mathrm{IA}(i)\setminus\mathrm{IA}(i+1)$ for some $i\geq 1$.  Then 
\bes
\phi=(x_1+f_1+F_1,\ldots,x_n+f_n+F_n), 
\ees
where $f_j\in A_{i+1}, F_i\in A_{i+2}+A_{i+3}+\ldots$ for all $1\leq j\leq n$. The derivation  
\bes
T(\phi)=f_1\partial_1+\ldots+f_n\partial_n\in L_i
\ees
will be called the {\em tangent to the automorphism $\phi$ with respect to power series topology}. Note that $T(\phi)\neq 0$.  Set also $T(\mathrm{Id})= 0$, where $\mathrm{Id}$ is the identity automorphism.

 For completness of the text we repeat some proofs given in \cite{Anick,BD93-2}.

\begin{lm}\label{l1} \cite{Anick,BD93-2} The following statements are true: 

(1) If $\phi,\psi \in  \mathrm{IA}(i)$ for some $i\geq 1$, then $T(\phi\psi)= T(\phi)+T(\psi)$ if $T(\phi)+T(\psi)\neq 0$ and $T(\phi^{-1})=-T(\phi)$;

(2) If $\phi\in  \mathrm{IA}(i)$ and $\psi\in  \mathrm{IA}(j)$ for some $i,j\geq 1$, then    $T([\phi,\psi])= [T(\phi),T(\psi)]\in L_{i+j}$ if $[T(\phi),T(\psi)]\neq 0$, where $[\phi,\psi]=\phi^{-1}\psi^{-1}\phi\psi$ is the group commutator. 
\end{lm}
\Proof 
 Let $\phi \in \mathrm{IA}(i)$ and $\psi \in \mathrm{IA}(j)$ for some $i,j\geq 1$,  
\bes
\phi=(x_1+f_1+F_1,\ldots,x_n+f_n+F_n),  T(\phi)=f_1\partial_1+\ldots+f_n\partial_n, 
\ees
where $F_k\in A_{i+2}+A_{i+3}+\ldots$ for all $1\leq k\leq n$, and 
\bes
\psi=(x_1+g_1+G_1,\ldots,x_n+g_n+G_n), T(\psi)=g_1\partial_1+\ldots+g_n\partial_n, 
\ees
where $G_k\in A_{j+2}+A_{j+3}+\ldots$ for all $1\leq k\leq n$. 

Notice that for any $f\in A_m$ with $m\geq 1$ we have 
\bee\label{f5}
\phi(f)=f(\phi(x_1),\ldots,\phi(x_n))=f(x_1+f_1+F_1,\ldots,x_n+f_n+F_n)\\
\nonumber =f+T(\phi)(f)+ F,  \ F\in A_{m+i+1}+A_{m+i+2}+\ldots. 
\eee

Applying (\ref{f5}), we immediately obtain  
\bee\label{f6}
\phi\psi(x_k)=\phi(x_k)+\phi(g_k)+\phi(G_k)= x_k+f_k+F_k+g_k+G_k+T(\phi)(g_k) +H_k,  
\eee
where $ H_k\in A_{i+j+2}+A_{i+j+3}+\ldots$ for all $1\leq k\leq n$. 

If $i=j$ and $T(\phi)+T(\psi)\neq 0$  then (\ref{f6}) implies that $ \phi\psi\in\mathrm{IA}(i)$ and $T(\phi\psi)= T(\phi)+T(\psi)$. Moreover, if $\psi=\phi^{-1}$ then (\ref{f6}) implies that 
\bes
f_k+g_k=0, F_k+G_k+T(\phi)(g_k)\in A_{i+j+2}+A_{i+j+3}+\ldots, 
\ees
i.e., $T(\phi^{-1})=-T(\phi)$ and 
\bes
\phi^{-1}=(x_1-f_1-F_1+T(\phi)(f_1)+F_1',\ldots,x_n-f_n-F_n+T(\phi)(f_n)+F_n'),  
\ees
where $F_k'\in A_{2i+2}+A_{2i+3}+\ldots$ for all $1\leq k\leq n$. 

Similarly, 
\bes
\psi^{-1}=(x_1-g_1-G_1+T(\psi)(g_1)+G_1',\ldots,x_n-g_n-G_n+T(\psi)(g_n)+G_n'),  
\ees
where $G_k'\in A_{2j+2}+A_{2j+3}+\ldots$ for all $1\leq k\leq n$. 
By (\ref{f6}), 
\bes
\phi\psi(x_k)=\psi( x_k)+f_k+F_k+T(\phi)(g_k) +H_k. 
\ees
Then, using (\ref{f5}), we obtain 
\bes
\psi^{-1}\phi\psi(x_k)= x_k+\psi^{-1}(f_k)+\psi^{-1}(F_k)+\psi^{-1}(T(\phi)(g_k)) +\psi^{-1}(H_k)\\
=x_k+f_k-T(\psi)(f_k)+F_k+T(\phi)(g_k)+H_k'=\phi(x_k)-T(\psi)(f_k)+T(\phi)(g_k)+H_k', 
\ees
where $H_k'\in A_{i+j+2}+A_{i+j+3}+\ldots$ for all $1\leq k\leq n$. Furthermore 
\bes
\phi^{-1}\psi^{-1}\phi\psi(x_k)= x_k-\phi^{-1}(T(\psi)(f_k))+\phi^{-1}(T(\phi)(g_k))+\phi^{-1}(H_k')\\
= x_k-T(\psi)(f_k)+T(\phi)(g_k)+R_k , 
\ees
where $R_k\in A_{i+j+2}+A_{i+j+3}+\ldots$ for all $1\leq k\leq n$. Recall that 
\bes
[T(\phi),T(\psi)]=\sum_{k=1}^n (T(\phi)(g_k)-T(\psi)(f_k))\partial_k\in L_{i+j}. 
\ees
Consequently, if $[T(\phi),T(\psi)]\neq 0$ then $T([\phi,\psi])=[T(\phi),T(\psi)]$. $\Box$

Let $H$ be an arbitrary subgroup of $\mathrm{Aut}(A)$ containing the subgroup of all affine automorphisms $\mathrm{Aff}_n(K)$ if $\mathfrak{M}$ is a unitary variety of algebras. If $\mathfrak{M}$ is not unitary then we assume that $H$ contains the subgroup of all linear automorphisms $GL_n(K)$. For convenience of notation, we set $\mathrm{Gr}_n=\mathrm{Aff}_n(K)$ if $A$ is unital, and $\mathrm{Gr}_n=GL_n(K)$ otherwise. Thus $\mathrm{Gr}_n\subseteq H$.

Set $H_0=H$ and $H_i=H\cap \mathrm{IA}(i)$ for all $i\geq 1$. Then 
\bes
H=H_0\subseteq H_1\subseteq H_2\subseteq \ldots. 
\ees
 Set $V_i=V_i(H)=\{0\}\cup \{T(\phi) | \phi\in H_i\setminus H_{i+1}\}$ for all $i\geq 1$.

\begin{lm}\label{l2}  \cite{Anick,BD93-2} 
The set $V_i=V_i(H)$ is a vector space for each $i\geq 1$. 
\end{lm}
\Proof Let $\phi \in H_i\setminus H_{i+1}$ for some $i\geq 1$. If $\psi \in H_i\setminus H_{i+1}$  and $T(\phi)+T(\psi)\neq 0$ then $T(\phi\psi)= T(\phi)+T(\psi)\in V_i$ by Lemma \ref{l1}(1). Moreover, 
$T(\phi^{-1})=-T(\phi)\in V_i$ by the same lemma. Consequently, $V_i$ is an abelian group. 

  Let $F$ be the subset of elements $a\in K$ such that $a T(\phi)\in V_i$. Then $0\in F$ since $0\in V_i$. If $a\neq 0$ then this condition means that there exists $\psi\in H_i$ such that $T(\psi)=aT(\phi)$. Suppose that $a,b\in F$ and let $\psi_1,\psi_2\in H_i$ be such that $T(\psi_1)= a T(\phi)$ and $T(\psi_2)= b T(\phi)$. If $a\neq b$ then $T(\psi_1\psi_2^{-1})=(a-b)T(\phi)\neq 0$ by Lemma \ref{l2}(1). Consequently, $a-b\in F$ and $F$ is a subgroup of the additive of $K$. Since $1\in F$ it follows that $\mathbb{Z}\subseteq F$.  For any $0\neq \a\in K$ we have 
\bes
T((\a \mathrm{Id})\phi(\a \mathrm{Id})^{-1})=\a^i T(\phi)\in V_i. 
\ees
Consequently, $F$ contains all $i$th powers of elements of $K$. For any rational number $r/s$ we have $r/s= (rs^{i-1}) (s^{-1})^i\in F$.  Every element $c\in K$ can be written as a rational linear combination of $c^i,(c+1)^i,\ldots,(c+i)^i$. Consequently, $F=K$. Therefore $V_i$ is a vector space over $K$. $\Box$

\begin{lm}\label{l3}  
$[L_0,V_i]\subseteq V_i$ for each $i\geq 1$. 
\end{lm}
\Proof The Lie algebra $L_0=\mathfrak{gl}(n)$ is generated by $x_1\partial_1+\ldots+x_n\d_n$ and all $x_i\d_j$, $1\leq i\neq j\leq n$.

If $s=x_1\partial_1+\ldots+x_n\d_n$ and $v\in V_i$ then $[s,v]=i v\in V_i$ by Lemma \ref{l2}. 

Let $s=x_2\partial_1$. Consider the automorphism $\a=(x_1+tx_2,x_2,\ldots,x_n)$. If $T(\phi)=f_1\partial_1+\ldots+f_n\partial_n$ then 
\bes
T(\a\phi\a^{-1})=\a T(\phi)\a^{-1}=\a(f_1-tf_2)\partial_1+\a(f_2)\partial_2+\ldots+\a(f_n)\partial_n. 
\ees

For any $f\in A_{i+1}$ we have 
\bes
\a(f)=f+ts(f)+\frac{1}{2!}t^2s^2(f)+\ldots+\frac{1}{(i+1)!}t^{i+1}s^{i+1}(f). 
\ees
Consequently, 
\bes
T(\a\phi\a^{-1})=T(\phi)+tD_1+t^2D_2+\ldots+t^{i+2}D_{i+2}, 
\ees
where $D_k\in L_i$ for all $1\leq k\leq i+2$, and 
\bes
D_1=(s(f_1)-f_2)\partial_1+s(f_2)\partial_2+\ldots+s(f_n)\partial_n= [s,T(\phi)]. 
\ees
Since $V_i$ is a vector space, an application of the Vandermonde determinant argument yields that $D_k\in V_i$  for all $1\leq k\leq i+2$. 
In particular, $D_1=[s,T(\phi)]\in V_i$. 

Similarly, $[x_i\partial_j,V_i]\subseteq V_i$ for all $i\neq j$.

If $\psi\in H_j$ for some $j\geq 1$ and $[T(\phi),T(\psi)]\neq 0$ then  $T([\phi,\psi])= [T(\phi),T(\psi)]\in V_{i+j}$ by Lemma \ref{l1}(2). $\Box$

\begin{lm}\label{l4} 
$[L_{-1},V_i]\subseteq V_{i-1}$ for each $i\geq 2$. 
\end{lm}
\Proof  Let $\phi\in H_i\setminus H_{i+1}$ where $i\geq 2$. It follows that  
\bes
\phi=X+F_{i+1}+F_{i+2}+\ldots+F_k, 
\ees
where $F_p$ is a homogeneous $n$-tuple of degree $p$ for all $s$ with $i+1\leq p\leq k$, and $T(\phi)=D_{F_{i+1}}\neq 0$. 

 Let $s=\partial_1\in L_{-1}$. Consider the automorphism $\a=(x_1+t,x_2,\ldots,x_n)$. We obtain  
\bes
\a\phi\a^{-1}=X+\sum_{p=i+1}^k\a(F_{p}). 
\ees
Then 
\bes
\a(F_{p})=F_{p}+ts(F_{p})+\frac{1}{2!}t^2s^2(F_{p})+\ldots+\frac{1}{p!}t^{p}s^{p}(F_{p})
\ees
and $\deg(s^r(F_{p}))=p-r$ if $r\leq p$. Consequently, 
\bes
\a\phi\a^{-1}=G_{0}+G_1+G_2+\ldots+G_{i+1}+\ldots G_k, 
\ees
where $G_k=F_k$,
\bes
 G_p=F_p+ts(F_{p+1})+\frac{1}{2!}t^2s^2(F_{p+2}) +\ldots+\frac{1}{(k-p)!}t^{k-p}s^{k-p}(F_{k}),  
\ees
where $i+1\leq p\leq k$, 
\bes
G_r=\frac{1}{(i+1-r)!}t^{i+1-r}s^{i+1-r}(F_{i+1})+\ldots+\frac{1}{(k-r)!}t^{k-r}s^{k-r}(F_k), 
\ees
where $0\leq r\leq i$ and $r\neq 1$, and 
\bes
G_1=X+\frac{1}{i!}t^{i}s^{i}(F_{i+1})+\ldots+\frac{1}{(k-1)!}t^{k-1}s^{k-1}(F_k). 
\ees
Notice that we can get rid of the constant part $G_0$ of $\a\phi\a^{-1}$ by considering the automorphism $\phi_1=\a\phi\a^{-1} (X-G_0)\in H$. We have 
\bes
\phi_1=G_1+G_2+\ldots+G_{i+1}+\ldots G_k, \ \ \ G_1=X+t^2G_1', 
\ees
and 
\bes
\phi_1=\phi+t(s(F_{i+1})+\ldots+s(F_{k}))+t^2\Phi_1. 
\ees
The inverse $G_1^{-1}$ of the linear automorphism $G_1$ is equal to the linear part of the automorphism $ (X+G_0)\a\phi^{-1}\a^{-1}$. Since $ (X+G_0)\a\phi^{-1}\a^{-1}\in \mathrm{Aut}_{K[t]}A[t]$ it follows that $G_1^{-1}\in \mathrm{Aut}_{K[t]}A[t]$. Obviously, $G_1^{-1}$ has the form 
\bes
G_1^{-1}=X+t^2G_1''
\ees
and
\bes
\phi_2=G_1^{-1}\phi_1=\phi+t(s(F_{i+1})+\ldots+s(F_{k}))+t^2\Phi_2. 
\ees
The linear part of $\phi_2$ is equal to $X$.  

More generally, suppose that for some $2\leq j\leq i$ we can find an automorphism $\phi_j\in H$ such that 
\bee\label{f9}
\phi_j=\phi+t^{r_j}(s(F_{i+1})+\ldots+s(F_{k}))+t^{r_j+1}\Phi_j  
\eee
and
\bee\label{f10}
\phi_j=X+E_j+\ldots+ E_k+\ldots, 
\eee
where $E_p$ is homogeneous of degree $p$ with respect to $x_1,\ldots,x_n$ for all $2\leq j\leq p$ and $E_j\neq 0$. 
Substituting $t^{j-1}$ instead of $t$ in (\ref{f9}) if necessary, we can assume that all powers of $t$ in (\ref{f9}) are divisable by $j-1$. 

Let $j<i$. Then  $t^{r_{j}+1}|E_j$ and let 
$E_j= t^{i_1}E_j^{(1)}+\ldots + t^{i_r}E_j^{(r)}$, where $r_{j}+1\leq i_1<\ldots <i_r$.  
Then 
\bes
T(\phi_j)=T_{E_j}=t^{i_1}T_{E_j^{(1)}}+\ldots + t^{i_r}T_{E_j^{(r)}}\in V_{j-1}(H). 
\ees
By Lemma \ref{l2}, $V_{j-1}(H)$ is a vector space over $K$. An application of the Vandermonde determinant argument yields that $T_{E_j^{(q)}}\in V_{j-1}(H)$  for all $1\leq q\leq r$. Consequently, there exists $\psi_q\in H$ such that 
\bes
\psi_q=X+E_j^{(q)} +C_{jq}, 
\ees
where the coordinates of the $n$-tuple $C_{jq}$ belong to $A_{j+1}+A_{j+2}+\ldots$. 

Let $i_q=(j-1)\mu_q$. Set $\a_q (x_i)=t^{\mu_q}x_i$ for all $i$. Then 
\bes
\theta_q=\a_q\psi_q\a_q^{-1}=X+t^{i_q}E_j^{(q)}+t^{i_q+1}C_{jq}'.
\ees
Set $\phi_{j+1}=\theta_1^{-1}\ldots \theta_r^{-1}\phi_j$. Notice that $\phi_{j+1}$ has the same form (\ref{f9}) since $r_{j}+1\leq i_1<\ldots <i_r$. By construction, we eliminated $E_j$ in (\ref{f10}) and $\phi_{j+1}$ can be written in the form  
\bes
\phi_{j+1}=X+E_{j+1}'+\ldots+ E_k'+\ldots, 
\ees
where $E_p'$ is homogeneous of degree $p$ with respect to $x_1,\ldots,x_n$ for all $j+1\leq p$. 

Suppose that $j=i$ in (\ref{f10}). Then 
\bes
E_j=t^{r_j}(s(F_{i+1}))+t^{r_j+1}E_j'. 
\ees
Repeating the same discussions as above together with an application of the Vandermonde determinant argument, we obtain $D_{s(F_{i+1})}=[s,D_{F_{i+1}}]=[s,T(\phi)]\in V_{i-1}(H)$. $\Box$

Set 
\bes
T(H)=L_{-1}\oplus L_0\oplus V_1\oplus V_2\oplus\ldots \oplus V_k\oplus\ldots. 
\ees
Recall that $L_{-1}=0$ if $A$ is not unital, that is, in this case  
\bes
T(H)=L_0\oplus V_1\oplus V_2\oplus\ldots \oplus V_k\oplus\ldots. 
\ees
\begin{theor}\label{t}
$T(H)$ is a graded Lie subalgebra of the Lie algebra $\mathrm{Der}(A)$.
\end{theor}
\Proof Lemmas (\ref{l2})--(\ref{l4}) directly imply the statement of the theorem. $\Box$

We call $T(H)$ the {\em tangent algebra} of $H$ with respect to power series topology. 

Notice that the groups of automorphisms of the polynomial algebra $K[x,y]$ and the free associative algebra $K\langle x,y\rangle$ are isomorphic. It is not difficult to check that the tangent algebras of $\mathrm{Aut}K[x,y]$ and $\mathrm{Aut}K\langle x,y\rangle$ are not isomorphic. The power series topologies of $\mathrm{Aut}K[x,y]$ and $\mathrm{Aut}K\langle x,y\rangle$ differ. 

The definition of $T(H)$ provided in this section is especially useful for establishing its connection to the study of automorphisms.

\section{Universal enveloping algebras, Jacobian matrices,  and divergence}

\hspace*{\parindent}

Let $B$ be an arbitrary algebra of $\mathfrak{M}$. The universal enveloping algebra $U(B)=U_{\mathfrak{M}}(B)$ is the associative algebra with identity generated by all universal operators of left multiplication $l_x$ and right multiplication $r_x$ for all $x\in B$ \cite{Jacobson}. 
 Every  $B$-bimodule $V$ in the variety of algebras $\mathfrak{M}$ can be regarded as a left $U$-module with respect to the action 
\bes
L_bm=bm, R_bm=mb, \ b\in B, m\in M. 
\ees
Conversely, every left $U$-module can be considered as a $B$-bimodule  in the variety of algebras $\mathfrak{M}$ \cite{Jacobson}. 

If $\mathfrak{M}$ is a unital variety and $B$ is an algebra with identity $1$, then $l_1=r_1=\mathrm{Id}=1$. 

We give a short explanation of the structure of $U=U(A)$ when $A=K_\mathfrak{M}\langle x_1,\ldots,x_n\rangle$ is a free algebra. Consider $A$ as the subalgebra of the free algebra $B=K_\mathfrak{M}\langle x_1,\ldots,x_n,y_1,\ldots,y_n\rangle$ generated by $x_1,\ldots,x_n$. Denote by $\Omega_{A}$ the set of all homogeneous elements of $B$ of degree $1$  with respect to the set of variables $y_1,\ldots,y_n$. For any $a\in A$ the operators 
\bes
L_a, R_a : \Omega_{A}\to \Omega_{A}, L_a(w)=aw, R_a(w)=wa, w\in \Omega_{A}, 
\ees
are well defined. The algebra generated by all $L_a, R_a$, where $a\in A$, is the algebra $U=U(A)$. 
Notice that 
\bes
\Omega_{A}=Uy_1\oplus \ldots \oplus Uy_n
\ees
is a free $U$-module with free generators $y_1,\ldots,y_n$. The linear map 
\bes
D : A \to \Omega_{A}
\ees
defined by $D(x_i)=y_i$ for all $i$ and such that 
\bes
D(ab)=aD(b)+D(a)b= L_a D(b)+R_b D(a)
\ees
for all $a,b\in A$, is called the {\em universal derivation} of $A$ and $\Omega_{A}$ is called the {\em universal differential module} of $A$ \cite{U94,U96}. Notice that over a field of charateristic zero $D(a)=0$ if and only if $a\in K$. 

For every $a\in A$ there exist unique elements $u_1,u_2,\ldots,u_n\in U$ such that
\bes
D(a)=u_1y_1+u_2y_2+\ldots+u_ny_n.
\ees
The elements $u_i=\frac{\partial a}{\partial x_i}$ are called the Fox derivatives of $a\in A$ \cite{U94,U96}. Set also 
\bes
\d(a)=(\frac{\partial a}{\partial x_1},\ldots,\frac{\partial a}{\partial x_n}). 
\ees
Let $Y=(y_1,\ldots,y_n)^t$ where $t$ is the transpose. 
Then 
\bes
D(a)=\d(a) Y. 
\ees
For any $f,g_1,\ldots,g_n\in A$ we get the chain rule  
\bes
D(f(g_1,\ldots,g_n))= \frac{\partial f}{\partial x_1}(g_1,\ldots,g_n)D(g_1)+\ldots+\frac{\partial f}{\partial x_n}(g_1,\ldots,g_n)D(g_n). 
\ees

Let  $\phi=(f_1,f_2,\ldots,f_n)$ be an arbitrary endomorphism of $A$. Denote by 
\bes
J(\phi)=[\partial_j(f_i)]_{1\leq i,j\leq n}
\ees
its Jacobian matrix. We have 
\bes
\begin{bmatrix}
D(f_1)\\
\cdot \\
\cdot \\
\cdot \\
D(f_n)
\end{bmatrix}
=J(\phi)Y. 
\ees
If $\psi=(g_1,\ldots,g_n)$ is another endomorphism then 
\bes
\phi\circ\psi=(g_1(f_1,\ldots,f_n),\ldots,g_n(f_1,\ldots,f_n)). 
\ees
The chain rule immediately implies 
\bee\label{fj}
J(\phi\circ\psi)=\phi(J(\psi))J(\phi). 
\eee
This gives that $J(\phi)$ is invertible over $U$ if $\phi$ is an automorphism. 
 An analogue of the Jacobian Conjecture for the free algebra $A$ can be formulated as follows: If $J(\phi)$ is invertible over $U$, then $\phi$ is an automorphism of $A$. This conjecture is true for free nonassociative algebras, for free commutative (characteristic $\neq 2$) and anticommutative algebras \cite{Yagzhev2} (see also \cite{U94,U96}), for free associative algebras \cite{Dicks,Schofield},  and for free Lie algebras and superalgebras \cite{Reutenauer,Shpilrain,U93,ZM}. 

Recall that the left-symmetric product is defined by $u\d_i\cdot v\d_j=(u\d_i)(v)\d_j$. 
We have $[D_1,D_2]=D_1\cdot D_2-D_2\cdot D_1$ for any $D_1,D_2\in D$ and 
for any two $n$-tuples $F$ and $G$ the product of $D_F$ and $D_G$ (see  (\ref{f4})) is given by 
\bes
D_F\cdot D_G=D_{D_F(G)}. 
\ees
Set also $J(D_F)=J(F)$. 

Every derivation $D$ of $A$ can be uniquely extended to a derivation $D^*$ of $U$ by $D^*(L_a)=L_{D(a)}$, $D^*(R_a)=R_{D(a)}$, for all $a\in A$.

\begin{lm}\label{l5} \cite{U16TG} Let $D_1,D_2\in \mathrm{Der}(A)$. Then the following statements are true:

(i) $J(D_1\cdot D_2)=D_1^*(J(D_2))+J(D_2)J(D_1)$;

(ii) $J([D_1,D_2])=D_1^*(J(D_2))-D_2^*(J(D_1))-[J(D_1),J(D_2)]$. 
\end{lm}
\Proof Let $D_1=D_F$ and $D_2=D_G$. Then $D_1\cdot D_2=D_{D_F(G)}$. Let $X=(x_1,\ldots,x_n)$ and consider the endomorphism $(X+tF)=(x_1+tf_1,\ldots,x_n+tf_n)$ where $t$ is an independent parameter. Obviously,
\bes
(X+tF)\circ G=G+tD_F(G)+t^2G_2+\ldots.
\ees
Consequently, $D_F(G)=\frac{\partial}{\partial t} ((X+tF)\circ G)|_{t=0}$.
By (\ref{fj}), we get
\bes
J((X+tF)\circ G) =(X+tF)(J(G))J(X+tF)\\
= (J(G)+tD_F(J(G))+t^2T_2+\ldots)(I+tJ(F))\\
= J(G)+t(D_F(J(G))+J(G)J(F))+t^2M_2+\ldots.
\ees
Hence
\bes
J(D_F(G))=\frac{\partial}{\partial t} J((X+tF)(G))|_{t=0}=D_F^*(J(G))+J(G)J(F),
\ees
which proves (i). This, in turn, directly implies (ii). $\Box$

The concept of the divergence of derivations, as it appears in the context of free Lie algebras \cite{ES11,Satoh12} and free associative algebras \cite{AKKN}, is closely connected to the study of $\mathrm{IA}$-automorphisms of free groups. A general notion of divergence for derivations of arbitrary free algebras was introduced by Powell \cite{Powell21}. In this work, we present a slight generalization of this notion, designed to facilitate the study of automorphisms in the general case. 

\begin{df} Let $D=D_F$ be a derivation of $A$ for some $F=(f_1,\ldots,f_n)$.  The  divergence of  $D$, denoted by $\mathrm{div}(D)$, is the image of 
\bes
\frac{\d f_1}{\d x_1}+\frac{\d f_2}{\d x_2}+\ldots+\frac{\d f_n}{\d x_n}=\mathrm{Tr}(J(D))=\mathrm{Tr}(J(F))
\ees
in the quotient space 
\bes
U/([U,U]+R), 
\ees
where $R=\mathrm{Rad}(U)$ is the Jacobson radical of $U=U(A)$. 
\end{df}

The image of the same element in the quotient space $U/[U,U]$ is referred to as the divergence in \cite{AKKN,ES11,Powell21,Satoh12}. However, this notion is not well suited for studying automorphisms of certain free algebras. Even in the case of free metabelian algebras \cite{BD93-2,U93,U95SMJ}, the Fox derivatives and Jacobian matrices are typically defined over $U/R$.

The following property of the divergence is established in \cite[Lemma 8]{U16TG} (see also \cite[Theorem 8.2]{Powell21})
\begin{co}\label{c1} Let $D_1,D_2\in \mathrm{Der}(A)$. Then 
\bes
\mathrm{div}([D_1,D_2])=D_1^*(\mathrm{div}(D_2))-D_2^*(\mathrm{div}(D_1))
\ees
\end{co}
\Proof Since $\mathrm{Tr}$ is linear, by Lemma \ref{l5} we obtain 
\bes
\mathrm{Tr}(J([D_1,D_2]))=D_1^*(\mathrm{Tr}(J(D_2)))-D_2^*(\mathrm{Tr}(J(D_1)))-\mathrm{Tr}[J(D_1),J(D_2)]. 
\ees
Obviously, $\mathrm{Tr}[J(D_1),J(D_2)]\in [U,U]$. $\Box$

This corollary implies that the space 
\bes
\mathrm{SDer}(A)=S_{-1}\oplus S_0\oplus S_1\oplus \ldots \oplus S_k\oplus \ldots 
\ees
of all derivations of $A$ with zero divergence forms a subalgebra of $L=\mathrm{Der}(A)$. Notice that $S_{-1}=L_{-1}$ and $S_0=\mathrm{sl}_n(K)$ is the subalgebra of traceless matrices of $L_0=\mathrm{gl}_n(K)$. 

This algebra will be called the {\em special Lie algebra of derivations} of $A$ as in the case of polynomial algebras. 

By Corollary \ref{c1}, the space 
\bes
\widetilde{\mathrm{SDer}}(A)=S_{-1}\oplus L_0\oplus S_1\oplus \ldots \oplus S_k\oplus \ldots 
\ees
of all derivations of $A$ with constant divergence also is a subalgebra of $L$.

If $A$ doesn't have the identity element then $S_{-1}=L_{-1}=0$ and we obtain  
\bes
\mathrm{SDer}(A)=S_0\oplus S_1\oplus \ldots \oplus S_k\oplus \ldots, 
\ees
and 
\bes
\widetilde{\mathrm{SDer}}(A)=L_0\oplus S_1\oplus \ldots \oplus S_k\oplus \ldots 
\ees

In order to relate these algebras with automorphism groups we need a minor generalization of Theorem 3.7 from \cite{BD93-1}. 

Let $F_i$ be free associative algebras in the variables $y_{ij}$ for all $1\leq i\leq k$, $j\geq 1$, such that $\deg(y_{ij})=r_{ij}\in \mathbb{Z}_{\geq 1}$. Set $Q=F_1\otimes F_2\otimes\ldots \otimes F_k$. 

If $M_{s_i}(K)$ is the matrix algebra of order $s_i$ over $K$ for all $1\leq i\leq k$ then 
\bes
M_{s_1}(K)\otimes M_{s_2}(K)\otimes\ldots \otimes M_{s_k}(K)\simeq M_{s_1\ldots s_k}(K). 
\ees

The following lemma is a minor generalization of \cite[Proposition 3.7]{BD93-1} and \cite[Lemma 3.4]{Liu}. 

\begin{lm}\label{l6} Let $f\in Q=F_1\otimes \cdots\otimes F_k$ be an element such for any $s_1,\ldots,s_k\geq 1$ and for any homomorphisms $\phi_i: F_i\to M_{s_i}(K)$ we have $\mathrm{Tr}(\theta(f))=0$ where 
\bes
\theta=\phi_1\otimes \cdots\otimes \phi_k: Q=F_1\otimes \cdots\otimes F_k\to M_{s_1}(K)\otimes\ldots \otimes M_{s_k}(K)\simeq M_{s_1\ldots s_k}(K). 
\ees
 Then $f\in [Q,Q]$.
\end{lm}
\Proof  Notice that an element $f_i\in F_i$ is called {\em balanced} in \cite{BD93-1} if it is a linear combination of elements of the form $uv-vu$, which is exactly equivalent to $f_i\in [F_i,F_i]$. 

Set $Q_1=F_2\otimes \cdots\otimes F_k$. Notice that 
\bes
[Q,Q]= [F_1,F_1]\otimes Q_1+F_1\otimes [Q_1,Q_1]. 
\ees
Consequently, if $f\notin [Q,Q]$ then 
\bes
f=u_1\otimes g_1+u_2\otimes g_2+\ldots+u_t\otimes g_t, 
\ees
where $g_1,\ldots,g_t\in Q_1\setminus [Q_1,Q_1]$ and $u_1,\ldots,u_t\in F_1$ are linearly independent modulo $[F_1,F_1]$. 
Leading an induction on $k$, we may assume that there exists a homomorphism 
\bes
\theta' : Q_1\to  M_{s_2}(K)\otimes\ldots \otimes M_{s_k}(K)
\ees
 such 
$\mathrm{Tr}(\theta'(g_1))=\a_1\neq 0$. Set $\mathrm{Tr}(\theta'(g_i))=\a_i$ for all $2\leq i\leq t$. 
Since $\mathrm{Tr}(A\otimes B)=\mathrm{Tr}(A) \mathrm{Tr}(B)$ it follows that 
\bes 
\a_1 \mathrm{Tr}(u_1)+\ldots+\a_t \mathrm{Tr}(u_t)=0
\ees
is a nontrivial trace identity of the matrix algebra $M_{s_1}$ for all $s_1\geq 1$. This contradicts to the well known Razmyslov-Prochesi Theorem \cite{Prochesi76,Razmyslov74},  which implies that every nontrivial trace identity of the matrix algebra $M_q(K)$ has degree $\geq q$. $\Box$

\begin{lm}\label{l7} Let $Q=F_1\otimes \cdots\otimes F_k$. Let $A\in GL_n(Q)$ be an invertible matrix of the form  
\bes
A=I_n+A_p+A_{p+1}+\cdots+A_q, 1\leq p\leq q, 
\ees 
where $A_i$ is a matrix of homogeneous elements of degree $i$ for all $p\leq i\leq q$. Then $\mathrm{Tr}(A_p)\in [Q,Q]$.
\end{lm}
\Proof Set $f=\mathrm{Tr}(A_p)$. Recall that $r_{ij}=\deg(y_{ij})$. Let $\phi_i : F_i\to M_{s_i}(K[t])$ be a homomorphism such that $\phi_i(y_{ij})=t^{r_{ij}}A_{ij}$ if $A_p$ depends on $y_{ij}$ and $\phi_i(y_{ij})=0$ otherwise, where $A_{ij}\in M_{s_i}(K)$. Consider 
\bes
\theta=\phi_1\otimes \cdots\otimes \phi_k: Q\to M_{s_1}(K[t])\otimes\ldots \otimes M_{s_k}(K[t]) \hookrightarrow M_{s_1\ldots s_k}(K[t]). 
\ees
Then 
\bes
\theta(A)=I_{ns_1\ldots s_k}+t^pA_p'+\ldots+t^qA_q'\in M_{ns_1\ldots s_k}(K[t])
\ees
and 
\bes
\mathrm{Tr}(\theta(f))=t^p\mathrm{Tr}(A_p'). 
\ees
Notice that $\det(\theta(A))=1$ since $\theta(A)$ is invertible. On the other hand, 
\bes
\det(\theta(A))=1+t^p\mathrm{Tr}(A_p')+t^{p+1}a, a\in K[t]. 
\ees
Consequently, $\mathrm{Tr}(A_p')=0$ and $\mathrm{Tr}(\theta(f))=0$. Substituting $t=1$, we obtain that $f$ satisfies the conditions of Lemma \ref{l6}. Consequently, $f\in [Q,Q]$. $\Box$

\begin{theor}\label{t1} Let $\mathfrak{M}$ be  one of the following varieties of algebras: 

 (1) A Nielsen-Schreier variety of algebras;  

(2) The variety of all associative algebras;

(3) the variety of all associative and commutative algebras; 

(4) the variety of all metabelian Lie algebras. 

If $H$ is a subgroup of $\mathrm{Aut}(A)$ containing $\mathrm{Gr}_n$, then the tangent algebra $T(H)$ is a subalgebra of the algebra $\widetilde{\mathrm{SDer}}(A)$ of derivations of $A$ with constant divergence. 
\end{theor}
\Proof 
 If $\mathfrak{M}$ is a Nielsen-Schreier variety of algebras then $U(A)$ is a free associative algebra  \cite{U94,U96}. 

In the case of free associative algebras, we have $U(A)\simeq A\otimes A^{\mathrm{op}}$, that is, a tensor product of free associative algebras. 

If $A$ is a polynomial algebra, then $U(A)\simeq A$ is a polynomial algebra. Note that $A$ is a tensor power of polynomial algebras in one variable, which are again free associative algebras. 

If $A$ is a free metabelian algebra, then the identity 
\bes
[[x,y],[z,t]]=0
\ees
implies 
\bes
L_{[x,y]}L_z=0. 
\ees
Consequently, the ideal $R$ of $U$ generated by all $L_{[x,y]}$ is nilpotent of index $2$. Obviously, $U/R$ is a polynomial algebra in $n$ variavles, and $R=\mathrm{Rad}(U)$.

Suppose that $\phi\in H_i$ for some $i\geq 1$, $\phi=X+F_{i+1}+\ldots$, and $T(\phi)=D_{F_{i+1}}\neq 0$. Then  
\bes
J(\phi)=I+J(F_{i+1})+\ldots. 
\ees
By Lemma \ref{l7}, $\mathrm{Tr}(J(F_{i+1}))\in [U,U]+R$. Consequently,  $\mathrm{div}(T(\phi))=\mathrm{div}(D_{F_{i+1}})=0$. This means that $T(\phi)\in \widetilde{\mathrm{SDer}}(A)$. $\Box$

Thus $T(\mathrm{Aut}(A))$ is a subalgebra of $\widetilde{\mathrm{SDer}}(A)$ for certain varieties of algebras. The radical was incorporated into the definition of divergence precisely to ensure this statement holds. However, whether this remains true in general is still an open question.
\begin{prob} 
Is it always true that $T(\mathrm{Aut}(A))$ is a subalgebra of $\widetilde{\mathrm{SDer}}(A)$? 
\end{prob}

\section{Approximately tame automorphisms and density}

\hspace*{\parindent}

Consider the descending central series (\ref{f2}). An automorphism $\phi\in \mathrm{Aut}(A)$ is called {\em approximately tame} if there exists a sequence of 
tame automorphisms $\{\psi_k\}_{k\geq 0}$ such that $\phi\psi_k^{-1}\in \mathrm{IA}(k)$. The topology defined by this definition on $\mathrm{Aut}(A)$ is called the {\em formal power series topology}  \cite{Anick}. An automorphism $\phi$ will be called {\em absolutely wild} if it is not approximately tame.

It is well known that (\cite[Theorem 5.2.1]{Essen},\cite{Yagzhev}) if $n\geq 3$ and $K$ is a field of the characteristic $0$ then the group of all tame automorphisms $\mathrm{TAut}(K[x_1,\ldots,x_n])$ of the polynomial algebra $K[x_1,\ldots,x_n]$ is generated by all affine automorphisms $\mathrm{Aff}_n(K)$ and by the quadratic automorphism 
\bes
\epsilon=(x_1+x_2^2,x_2,\ldots,x_n). 
\ees
Moreover, Anick \cite{Anick} proved that the same subgroup is dense in the group of all automorphisms 
$\mathrm{Aut}(K[x_1,x_2,\ldots,x_n])$. Consequently, every automorphism of $K[x_1,x_2,\ldots,x_n]$ is approximately tame. 

Notice that  $\mathrm{Aut}(K[x_1,x_2])$ is not even finitely generated modulo $\mathrm{Aff}_2(K)$ (see, for example \cite{Bodnarchuk}).  Bodnarchuk \cite{Bodnarchuk} proved that  $\mathrm{TAut}(K[x_1,\ldots,x_n])$ for $n\geq 3$  is generated by all affine automorphisms $\mathrm{Aff}_n(K)$ and by any non-affine triangular automorphism.

\begin{lm}\label{l8} Let $H$ be an arbitrary subgroup of $\mathrm{Aut}(A)$ containing the subgroup  $\mathrm{Gr}_n$. Let $\phi\in \mathrm{IA}(i)\setminus\mathrm{IA}(i+1)$ for some $i\geq 1$. If $T(\phi)\in T(H)$ then there exists an automorphism $\psi\in H$ such that $\phi\psi^{-1}\in \mathrm{IA}(i+1)$.
\end{lm}
\Proof We have $T(\phi)\in V_i$. By the definition of $V_i$ there exists $\psi\in H_i\setminus H_{i+1}$ such that $T(\phi)=T(\psi)$. By Lemma \ref{l1}(1), $\phi\psi^{-1}\in \mathrm{IA}(i+1)$. $\Box$

More generally, let $H$ and $N$ be subgroups of $\mathrm{Aut}(A)$ such that $\mathrm{Gr}_n\subseteq H\subseteq N$. We say that {\em $H$ is {\em dense} in $N$} with respect to the power series topology if for any $i\geq 1$ and for any $\psi\in N\cap \mathrm{IA}(i)$ there exists $\phi\in H$ such that $\phi^{-1}\psi\in N\cap \mathrm{IA}(i+1)$.

\begin{co}\label{c2} Let $\mathrm{Gr}_n\subseteq H\subseteq N$ be subgroups of $\mathrm{Aut}(A)$.  If $T(H)=T(N)$ then $H$ is dense in $N$. 
\end{co}
\Proof Recall that $N_i=N\cap \mathrm{IA}(i)$ and $V_i(N)=\{0\}\cup \{T(\phi) | \phi\in N_i\setminus N_{i+1}\}$ for all $i\geq 1$. Let $ \psi\in N_i\setminus N_{i+1}$. Then 
$0\neq T(\psi)\in V_i(N)$. Since $T(H)=T(N)$ it follows that $V_i(H)=V_i(N)$ and there exists $\phi\in H_i\setminus H_{i+1}$ such that  $T(\phi)=T(\psi)$. Consequently, 
$\phi\psi^{-1}\in V_{i+1}(N)$. 
$\Box$

\begin{theor}\label{t0}  If the tangent algebra $T(\mathrm{Aut}(A))$ is generated modulo $L_{-1}+L_0$ by all  derivations of the form $f\d_1$, where $f\in K_{\mathfrak M}\langle x_2,\ldots,x_n\rangle$ homogeneous of degree $\geq 2$, then every automorphism of $A$ is approximately tame. 
\end{theor}
\Proof Let $H=\mathrm{TAut}(A)$ be the group of all tame automorphisms of $A$. If $f\in K_{\mathfrak M}\langle x_2,\ldots,x_n\rangle$ is homogeneous of degree $\geq 2$, then 
\bes
f \d_1=T(x_1+f,x_2,\ldots,x_n)\in T(H). 
\ees
Consequently, $T(H)=T(\mathrm{Aut}(A))$. By Corollary \ref{c2}, $H$ is dense in $\mathrm{Aut}(A)$. 
$\Box$

This theorem raises the question of determining the generators of $T(\mathrm{Aut}(A))$. By Theorem \ref{t1}, $T(\mathrm{Aut}(A))$ is a subalgebra of $\widetilde{\mathrm{SDer}}(A)$ for certain varieties of algebras.

\begin{co}\label{c3} \cite{Anick} 
The subgroup $H$ of $\mathrm{Aut}(K[x_1,\ldots,x_n])$ generated by 
all affine automorphisms $\mathrm{Aff}_n(K)$ and by the quadratic automorphism $\epsilon$ is dense in $\mathrm{Aut}(K[x_1,\ldots,x_n])$. In particular, the group of all tame automorphisms is dense in the group of all automorphisms. 
\end{co}
\Proof Let $A=K[x_1,\ldots,x_n]$. It is well known (see, for example \cite{Shafarevich81}) that the subalgebra 
\bes
\widetilde{\mathrm{SDer}}(A)^{\geq 0}=L_0+S_1+S_2+\ldots
\ees
 of $\widetilde{\mathrm{SDer}}(A)$ modulo $L_0$ is generated by $T(\epsilon)=x_2^2\partial_1$. Consequently, $T(H)=\widetilde{\mathrm{SDer}}(A)=T(\mathrm{Aut}(A))$.   By Corollary \ref{c2},  $H$ is dense in $\mathrm{Aut}(A)$. $\Box$

One of the classical and well studied varieties of algebras is the variety of metabelian Lie algebras. Let $L_n$ be a free Lie algebra of rank $n$ in the variables $x_1,\ldots,x_n$ over a field $K$. 
Then $M_n = L_n/L''_n=L_n/[[L_n,L_n],[L_n,L_n]]$
 is the free metabelian Lie algebra 
of rank $n$ in the variables $y_i=x_i+L''_n$, where $1\leq i\leq n$. Obviously, any nontrivial exponential automorphism of $M_2$ is wild (see, for example \cite[Proposition 4]{Shmelkin73},\cite[Theorem 3]{Artamonov}, \cite[Theorem 3]{Papistas}) since every tame automorphism in this case is linear. For the same reason they are absolutely wild. In 1992, V. Drensky \cite{Drensky92} proved that the exponential automorphism 
\bes
\mathrm{exp}(\mathrm{ad}[y_1,y_2])=(y_1+[[y_1,y_2],y_1], y_2+[[y_1,y_2],y_2], y_3+[[y_1,y_2],y_3])
\ees
 of $M_3$ is wild.  In fact, the proof given in  \cite{Drensky92} shows that this automorphism is absolutely wild. 
More examples of non-tame automorphisms of $M_3$ are given in \cite{Nauryzbaev09,Romankov08}. 

 In 1993 Bryant and Drensky \cite{BD93-2} proved that every automorphism of $M_n$ for $n\geq 4$ is approximately tame, that is, the group of tame automorphisms is dense in $\mathrm{Aut}(M_n)$. More detailed versions of the results from \cite{BD93-2} are provided in \cite{KP16} by Kofinas and Papistas. Here, we present direct formulations of the results established in \cite{BD93-2,KP16}, expressed in the language of tangent algebras. 
\begin{theor}\label{t2} \cite{BD93-2,KP16} Let $M_n$ be a free Lie algebra over a field $K$ of characteristic zero in the variables $y_1,\ldots,y_n$. Then the following statements are true. 

(a) If $n\geq 4$ then the algebra $\widetilde{\mathrm{SDer}}(M_n)$ is generated modulo $L_0$ by the derivation $[y_2,y_3]\d_1$.

(b) The algebra $\widetilde{\mathrm{SDer}}(M_3)$ is generated modulo $L_0$ by the derivations $[y_2,y_3]\d_1$ and 
\bes
\mathrm{ad}([y_1,y_2])=[[y_1,y_2],y_1]\partial_1+[[y_1,y_2],y_2]\partial_2+[[y_1,y_2],y_3]\partial_3. 
\ees
\end{theor}
\Proof  Let $R$ be the subalgebra of $\widetilde{\mathrm{SDer}}(M_n)$ generated by $L_0$ and $[y_2,y_3]\d_1$. By Proposition \ref{X}, $R$ is homogeneous and 
\bes
R=L_0\oplus R_1\oplus R_2\oplus\ldots \oplus R_k\oplus\ldots. 
\ees
Furthermore, Proposition \ref{X} also shows that each $R_k$ is invariant under the action of the linear automorphism group $GL_n(K)$.

In fact, the main result of \cite{BD93-2,KP16} states that if $n\geq 4$, then the subalgebra of $\widetilde{\mathrm{SDer}}(M_n)^{\geq 1}=S_1\oplus S_2\oplus\ldots$, generated by 
$[y_2,y_3]\d_1$ and closed under the action of $GL_n(K)$, is equal to $\widetilde{\mathrm{SDer}}(M_n)^{\geq 1}$. Consequently, 
\bes
R_1\oplus R_2\oplus\ldots=S_1\oplus S_2\oplus\ldots.
\ees
 Therefore $R=\mathrm{SDer}(M_n)$. This proves $(a)$. 

Obviously, the same discussion applies to case $(b)$. $\Box$

There are many similarities and intersections between the results obtained in \cite{BD93-1} and those in  \cite{ES11}. The free metabelian Lie ring, commonly used in the literature, is referred to as the Chen ring in \cite{ES11}. 

It is not difficult to show that $\widetilde{\mathrm{SDer}}(M_2)$ is not finitely generated modulo $L_0$. 

\begin{co}\label{c4} \cite{BD93-2,KP16} (a) If $n\geq 4$ then the group of automorphisms generated by all linear automorphisms $GL_n(K)$  and  by the quadratic automorphism 
\bes
\tau=(y_1+[y_2,y_3],y_2,\ldots,y_n)
\ees
is dense in  $\mathrm{Aut}(M_n)$. In particular, the group of all tame automorphisms is dense in the group of all automorphisms. 

(b) If $n=3$ then the group of automorphisms generated by all linear automorphisms $GL_3(K)$, $\tau$,  and  $\mathrm{exp}(\mathrm{ad}[y_1,y_2])$ is dense in $\mathrm{Aut}(M_3)$. 
\end{co}
\Proof We have 
\bes
T(\tau)=[y_2,y_3]\d_1, \ \ T(\mathrm{exp}(\mathrm{ad}[y_1,y_2]))=\mathrm{ad}([y_1,y_2]). 
\ees
It remains to apply Corollary \ref{c2} and Theorem \ref{t2}. $\Box$

Recall that in 2015, Nauryzbaev \cite{Nauryzbaev15} proved that for $n\geq 4$ over a field of characteristic different from $3$, the group of tame automorphisms $\mathrm{TAut}(M_n)$  is generated by all linear automorphisms $GL_n(K)$ together with the quadratic automorphism $\tau$. It was recently shown in \cite{U24} that for
$n\geq 4$ over a field of characteristic different from $3$, the automorphism group $\mathrm{Aut}(M_n)$ is generated by all linear automorphisms  $GL_n(K)$, 
 the quadratic automorphism $\tau$, and the cubic Chein automorphism 
\bes
(y_1+[[y_2,y_3],y_1],y_2,\ldots,y_n). 
\ees
This is a close analogue of the well-known Bachmuth-Mochizuki-Roman'kov Theorem \cite{BM85,Romankov85} which states that every automorphism of a free metabelian group of rank $\geq 4$ is tame.

\section{Absolutely wild automorphisms and their detection}

\hspace*{\parindent}

 Obviously, every absolutely wild automorphism is wild, but not vice versa. For example, the Nagata automorphism is wild but not absolutely wild.  At the moment we don't know if the Anick automorphism is approximately tame or not. But almost all known to us other examples of wild automorphisms are absolutely wild. Here we demonstrate a method establishing of absolutely wild automorphisms, which is a clarification of methods used by V. Drensky \cite{Drensky92}, Papistas 
\cite{Papistas93}, and Bahturin and Shpilrain \cite{BSh95}. 

First, extend some notations for endomorphisms. Let $E=\mathrm{End}(A)$ be the monoid of all endomorphisms of $A$.  
Let $\mathrm{IE}(i)=\mathrm{IE}(A,i)$ be the submonoid of all endomorphisms of $A$ that induces the identity automorphism on the factor-algebra $A/(A_{i+1}+A_{i+2}+\ldots)$. 
Let $\phi\in \mathrm{IE}(i)\setminus\mathrm{IE}(i+1)$ for some $i\geq 1$.  Then 
\bes
\phi=(x_1+f_1+F_1,\ldots,x_n+f_n+F_n), 
\ees
where $f_j\in A_{i+1}, F_i\in A_{i+2}+A_{i+3}+\ldots$ for all $1\leq j\leq n$. Set   
\bes
T(\phi)=f_1\partial_1+\ldots+f_n\partial_n\in L_i
\ees
as in the case automorphisms.  Set also $T(\mathrm{id})= 0$. 

\begin{theor}\label{Y} Let $\mathfrak{M}$ be  one of the following varieties of algebras: 

 (1) A Nielsen-Schreier variety of algebras;  

(2) The variety of all associative algebras;

(3) the variety of all associative and commutative algebras; 

(4) the variety of all metabelian Lie algebras. 

Let $\mathfrak{N}$ be any subvariety of $\mathfrak{M}$ and let $I\subseteq A$ be the ideal of all identities of $\mathfrak{N}$ in $\mathfrak{M}$ in $n$ variables. Suppose that an endomorphism $\epsilon\in \mathrm{IE}_i(A)\setminus\mathrm{IE}_{i+1}(A)$ for some $i\geq 1$ induces an automorphism $\phi$ of $B=A/I$ and the ideal $I$ does not contain elements of degree $\leq i+1$.
 If $\mathrm{div}(T(\epsilon))\neq 0$ then $\phi$ is absolutely wild. 
\end{theor}
\Proof Let 
\bes
T(\epsilon)=f_1\partial_1+\ldots+f_n\partial_n\in L_i.
\ees
 Denote by $a'$ the image of any element $a\in A$ in $B$ under the natural projection. Then 
\bes
T(\phi)=f_1'\partial_1+\ldots+f_n'\partial_n \neq 0
\ees
since $I$ does not contain elements of degree $\leq i+1$ and 
\bes
\phi=(z_1+f_1'+F_1',\ldots,z_n+f_n'+F_n'), f_j\in A_{i+1}, F_j\in A_{i+2}+A_{i+3}+\ldots. 
\ees

 If $\phi$ is approximately tame then there exists a tame automorphism $\psi$ of $B$ such that $\phi\psi^{-1}\in \mathrm{IA}(B,i+1)$. Obviously, this implies that $\psi\in \mathrm{IA}(B,i)$ and $T(\phi)=T(\psi)$. Then 
\bes
\psi=(z_1+f_1'+G_1',\ldots,z_n+f_n'+G_n'),  G_j\in A_{i+2}+A_{i+3}+\ldots, 
\ees
Suppose that $\psi$ is induced by an automorphism $\tau$ of $A$. Then 
\bes
\tau=(x_1+f_1+G_1+T_1,\ldots,x_n+f_n+G_n+T_n),   \ T_j\in I, \ 1\leq j\leq n. 
\ees
 We have $\deg(T_j)\geq i+2$ since $I$ does not contain elements of degree $\leq i+1$.  Then 
\bes
\tau=(x_1+f_1+H_1,\ldots,x_n+f_n+H_n),  H_j\in A_{i+2}+A_{i+3}+\ldots.  
\ees
Therefore, $T(\tau)=T(\epsilon)$ and $\mathrm{div}\,T(\tau)=\mathrm{div}\,T(\epsilon)\neq 0$. This contradicts Theorem \ref{t1}. $\Box$

Let $\mathfrak{N}_{c}$ be the variety of all nilpotent Lie algebras of class $\leq c+1$, where $c\geq 1$, and let 
$\mathfrak{A}$ be the variety of all abelian Lie algebras. We have $\mathfrak{N}_1=\mathfrak{A}$. 
A variety of Lie algebras of the form $\mathfrak{N}_{c_1}\mathfrak{N}_{c_2}\ldots\mathfrak{N}_{c_k}$ is called polynilpotent.  In paticular, $\mathfrak{A}^2$ is the variety of all metabelian Lie algebras.In 1992, V. Drensky \cite{Drensky92}, in 1993, Papistas \cite{Papistas93}, and in 1995, Bahturin and Shpilrain \cite{BSh95} proved that free algebras of rank $n \geq 2$ in any polynilpotent variety $\mathfrak{M}$ of Lie algebras that is not $\mathfrak{A}$ and  $\mathfrak{A}^2$ have wild automorphisms. However, it should be noted that by “polynilpotent,” they in fact considered much broader classes of varieties. Their proofs can be used to show that those wild automorphisms are absolutely wild.  To clearly illustrate the method, we provide the following specific result together with its proof.

\begin{theor}\label{t3}
Let $\mathfrak{M}$ be a polynilpotent variety of Lie algebras that is not $\mathfrak{A}$ and $\mathfrak{A}^2$. 
Then every free algebra of $\mathfrak{M}$ of rank $n\geq 2$ has absolutely wild automorphisms. 
\end{theor}
\Proof Let $\mathfrak{M}=\mathfrak{N}_{c_k}\mathfrak{N}_{c_{k-1}}\ldots\mathfrak{N}_{c_1}$. If $k=1$ then $\mathfrak{M}=\mathfrak{N}_{c_1}$ and $c_1\geq 2$ since $\mathfrak{M}\neq \mathfrak{A}$. In this case the automorphism 
\bes
(x_1+[x_1,x_2],x_2,\ldots,x_n) 
\ees
is absolutely wild by Theorem \ref{Y}. 

Suppose that $k\geq 2$. Let $I$ be the ideal of identities of $\mathfrak{M}$ in the free Lie algebra $L_n$, i.e., all identities of $\mathfrak{M}$ in $n$ variables. 
Set  
\bes
J=(\ldots((L_n^{c_1+1})^{c_2+1})^{\ldots})^{c_{k-1}+1}. 
\ees 
Then 
\bes
I=(\ldots((L_n^{c_1+1})^{c_2+1})^{\ldots})^{c_k+1}=J^{c_k+1}
\ees
and $L_n/I=L_n/(J^{c_k+1})$ is the free algebra of $\mathfrak{M}$ of rank $n$. 

Set also
\bes 
u_1=\mathrm{ad}(x_1)^{c_1}(x_2)\in L_n^{c_1+1}, \ \ u_2=\mathrm{ad}(u_1)^{c_2}(\mathrm{ad}{x_1}(u_1))\in (L_n^{c_1+1})^{c_2+1}.
\ees
 Suppose that $u_t\in (\ldots((L_n^{c_1+1})^{c_2+1})^{\ldots})^{c_t+1}$ is already constructed for some $2\leq t<k$. Then set \bes
u_{t+1}=\mathrm{ad}(u_t)^{c_{t+1}}(\mathrm{ad}{x_1}(u_t))\in ((L_n^{c_1+1})^{c_2+1})^{\ldots})^{c_{t+1}+1}. 
\ees

Consider the lexocographic order on the set of basis elements of the free associative algebra $U(L_n)$ defined by $x_1>x_2>\ldots>x_n$. Denote by $\overline{a}$ the leading monomial of any homogeneous element $a\in U(L)$. 
Notice that 
\bes
\overline{u_1}=x_1^{c_1}x_2, \overline{u_{t+1}}=x_1\overline{u_t}^{c_{t+1}+1}, 
\ees
for all $t\geq 1$. Consequently,  
\bes
\deg(u_1)=c_1+1, \deg(u_{t+1})=\deg(u_t)(c_{t+1}+1)+1. 
\ees

Set $u=u_{k-1}\in J$. We show that 
\bee\label{99}
\deg(u)+2< (c_1+1)\ldots(c_k+1). 
\eee

If $k=2$ then $\deg(u)=c_1+1$ and (\ref{99}) does not hold only if $c_1=c_2=1$. This is impossible since $\mathfrak{M}\neq \mathfrak{A}^2$.

If $k\geq 3$ then 
\bes
\deg(u_2)+2=(c_1+1)(c_2+1)+3<(c_1+1)(c_2+1)(c_k+1)
\ees
By induction on $t$ we prove that 
\bes
\deg(u_t)+2<(\Pi_{i=1}^t(c_i+1))(c_k+1)
\ees
for all $2\leq t\leq k-1$. 
If it is true for some $t<k-1$ then 
\bes
\deg(u_t)(c_{t+1}+1)+2(c_{t+1}+1)<(\Pi_{i=1}^{t+1}(c_i+1))(c_k+1). 
\ees
We have 
\bes
\deg(u_{t+1})+2=\deg(u_t)(c_{t+1}+1)+3<\deg(u_t)(c_{t+1}+1)+2(c_{t+1}+1). 
\ees
Consequently, 
\bes
\deg(u_{t+1})+2 <(\Pi_{i=1}^{t+1}(c_i+1))(c_k+1). 
\ees
This proves  (\ref{99}). 

Let $u'$ be the image of $u$ in $L_n/I$. Then $\mathrm{ad}(u')$ is nilpotent and $\mathrm{exp}(\mathrm{ad}(u'))$ is a nonlinear automorphism of $L_n/I$ since $0\neq [u,x_1]\notin I$ by (\ref{f1}). 
If $n=2$ then $\mathrm{exp}(\mathrm{ad}(u'))$ is an absolutely wild automorphism since the group of tame automorphisms coincides with $GL_2(K)$.

Assume that $n\geq 3$. Let $w=u(x_2,x_3)$. Consider the endomorphism 
\bes
\psi = (x_1+[[w,x_1]x_1],x_2,\ldots,x_n) 
\ees
of $L_n$. 

First check that $\psi$ induces an automorphism $\phi$ of $L_n/I$. Obviously, $\psi$ induces an automorphism of $L_n/J^2$. Then  it induces an automorphism of $J^t/J^{t+1}$ for all $t\geq 1$. Consequently, $\psi$ induces an automorphism of $L_n/I=L_n/(J^{c_k+1})$ (see more details in \cite{U93}). 

Notice that $T(\phi)$ is induced by $D=[[w,x_1]x_1]\partial_1\in L_i$, where $i=\deg(u)+1$. The ideal $I$ does not contain elements of degree $\leq i+1$ since 
\bes
i+1=\deg(u)+2<\Pi_{i=1}^{k}(c_i+1)
\ees
by (\ref{f1}). We also have 
\bes
\frac{\partial}{\partial x_1}([[w,x_1]x_1])=[w,x_1]-x_1w, \ \ \mathrm{div}(D)\neq 0. 
\ees
By Theorem \ref{Y}, $\phi$ is an absolutely wild automorphisms of $L_n/I$. $\Box$

\begin{lm}\label{l9} Let $\phi\in \mathrm{IA}_i(A)\setminus\mathrm{IA}_{i+1}(A)$ for some $i\geq 1$ be an automorphism of the free associative algebra $A=K\langle x_1,x_2\rangle$. Then $T(\phi)([x_1,x_2])=0$. 
\end{lm}
\Proof Obviously, every elementary automorphism preserves the commutator $[x_1,x_2]$ up to a non-zero scalar. Since every automorphism of $A$ is tame \cite{Czer, ML70}, every automorphism preserves the commutator $[x_1,x_2]$ up to a non-zero scalar.
Consequently, we get $\phi([x_1,x_2])=[x_1,x_2]$ since $\phi\in \mathrm{IA}_i(A)$, i.e., 
\bes
[x_1,x_2]=[x_1,x_2]+T(\phi)([x_1,x_2])+F, F\in A_{i+3}+A_{i+4}+\ldots. 
\ees
Hence $T(\phi)([x_1,x_2])=0$. $\Box$

G. Bergman \cite{Bergman} proved that the endomorphism 
\bes
\b=(x_1+[x_1,x_2]^2, x_2)
\ees
of $A=K\langle x_1,x_2\rangle$ induces a wild automorphism of the free algebra of rank two of the variety of algebras $\mathrm{Var}(M_2(K))$ generated by the  matrix algebra $M_2(K)$. 

\begin{lm}\label{l10} The endomorphism $\b$ of $A=K\langle x_1,x_2\rangle$ induces an absolutely wild automorphism of the free algebra of rank two of $\mathrm{Var}(M_2(K))$. 
\end{lm}
\Proof Let $I$ be the ideal of all identities of $M_2(K)$ in two variables and $A/I$ is the free algebra of rank two in $\mathrm{Var}(M_2(K))$. The $T$-ideal of identities of $M_2(K)$ is generated \cite{Drensky81} by 
\bes
\sum_{\s\in S_4} x_{\s(1)}x_{\s(2)}x_{\s(3)}x_{\s(4)}, [[x_1,x_2]^2,x_3], 
\ees
 and it does not contain elements of degree $\leq 4$ in two variables. Suppose that $\b$ induces an approximately tame automorphism $\overline{\b}$ of $A/I$. Then there exists a tame automorphism $\a$ of $A$ such that $\overline{\a}^{-1}\overline{\b}\in \mathrm{IA}(A/I,4)$. Therefore,  
\bes
\a=(x_1+[x_1,x_2]^2 +F_1, x_2+F_2), F_i\in A_5+A_6+\ldots. 
\ees
This gives
\bes
T(\a)=[x_1,x_2]^2\partial_1. 
\ees
Then we have 
\bes
T(\a)([x_1,x_2])=[[x_1,x_2]^2,x_2]\neq 0,
\ees
which contradicts Lemma \ref{l9}.  $\Box$

 More examples of wild automorphisms of free algebras in two variables are given in \cite{Drensky90}.

\section*{Acknowledgments}
The first author is supported by FAPESP  grant 2024/14914-9 and CNPq grant 305196/2024-3.
The second author would like to thank the Max Planck Institute f\"ur Mathematik  for
its hospitality and excellent working conditions, where some part of this work has been done. His research is also supported by grant AP23486782 from the Ministry of Science and Higher Education of the Republic of Kazakhstan, as well as by FAPESP grant 2025/01683-1.

\end{document}